\documentclass[12pt]{amsart}
\usepackage{graphicx,amssymb,stmaryrd,amsfonts,epsfig,amsthm,a4,amsmath,mathrsfs}
\usepackage{vmargin,url}
\usepackage[latin1]{inputenc}
\psfigdriver{dvips}

\newtheorem{thm}{Theorem}
\newtheorem{cor}[thm]{Corollary}
\newtheorem{lem}[thm]{Lemma}
\newtheorem{prop}[thm]{Proposition}
\theoremstyle{definition}

\newtheorem{que}[thm]{Question}
\theoremstyle{remark}

\numberwithin{equation}{section}
\begin{document}

\title{A sofic group away from amenable groups}
\author{Yves Cornulier}%
\date{Nov.~19, 2009}%June 18, 2009}

%20E26  	Residual properties and generalizations; residually finite groups
%43A07  	Means on groups, semigroups, etc.; amenable groups

\begin{abstract}
We give an example of a sofic group, which is not a limit of amenable groups.
\end{abstract}

\maketitle

\section{Introduction}

Let $G$ be a finitely generated group. A sequence $(G_n)$ of finitely generated groups {\it converges} to $G$ \cite{Chb,Gri,Cha,CG} if there exists a finitely generated free group $F$, and normal subgroups $N$, $(N_n)$ of $G$ such that

\begin{itemize}
\item $F/N\simeq G$; $F/N_n\simeq G_n$ for all $n$;
\item $(N_n)$ converges to $N$, i.e. for all $x\in N$, resp. $y\in F-N$; eventually $x\in N_n$, resp. $y\notin N_n$.
\end{itemize}

The finitely generated group $G\simeq F/N$ is {\it isolated} if whenever such a situation occurs, eventually $N_n=N$.

{\it Sofic groups} were introduced by Gromov as ``groups whose Cayley graph is initially subamenable" \cite[p.157]{Gro} and by B. Weiss in \cite{Wei}.
It will be enough for us to know that the class of sofic groups is not empty and satisfies the following properties

\begin{itemize}
\item \textbf{(subgroups)} Subgroups of sofic groups are sofic;
%\item Amenable implies sofic;
\item \textbf{(direct limits)} A group is sofic if (and only if) all its finitely generated subgroups are sofic
\item \textbf{(amenable extensions)} If a group is sofic-by-amenable, i.e.~has an amenable quotient with sofic kernel, then it is sofic as well (in particular, amenable implies sofic).
\item \textbf{(marked limits)} If a finitely generated group $G$ is a limit of a sequence of sofic groups $(G_n)$, then $G$ is sofic as well. In particular, residually finite groups are sofic.
\end{itemize}

No group is known to fail to be sofic.
The following question was asked by Gromov \cite{Gro} and is also addressed by A. Thom \cite{Th} and V. Pestov \cite{Pe}.

\begin{que}
Is every finitely generated sofic group a limit of amenable groups (``initially subamenable")?
\end{que}

Examples of sofic groups that are not residually amenable were obtained in \cite{ES}, but by construction, these groups are limits of finite groups.
Gromov expected a negative answer to the general question and we confirm this expectation.

\begin{thm}
There exists a finitely presented, non-amenable, isolated, (locally residually finite)-by-abelian finitely generated group.\label{main}
\end{thm}

\begin{cor}
There exists a finitely presented sofic group that is not a limit of amenable groups.
\end{cor}

Of course, being non-amenable and isolated, a group cannot be limit of amenable groups. Besides, by the remarks above, a (locally residually finite)-by-abelian group has to be sofic (a group is locally $\mathscr{X}$ means that every finitely generated subgroup is $\mathscr{X}$). Therefore the corollary follows from the theorem. 

Examples of isolated groups were provided in \cite{CGP}. In that paper, examples of infinite isolated groups with Kazhdan's Property T were given \cite[Paragraph 5.4]{CGP}, but it is not known if they are sofic. Also, some lattices in non-linear semisimple groups with finite center, are known to be isolated \cite[Paragraph 5.8]{CGP}, but they appear as natural candidates to be examples of non-sofic groups (see the discussion in Section \ref{es}).

Fix a prime $p$, and let $\Gamma$ be the group of matrices

$$\begin{pmatrix}
  a & b & u_{02} & u_{03} & u_{04} \\
  c & d & u_{12} & u_{13} & u_{14} \\
  0 & 0 & p^{n_2} & u_{23} & u_{24} \\
  0 & 0 & 0 & p^{n_3} & u_{34} \\
  0 & 0 & 0 & 0 & 1 \\
\end{pmatrix},$$

with 
$$\begin{pmatrix}
  a & b \\
  c & d \\
\end{pmatrix}\in\textnormal{SL}_2(\mathbf{Z}),u_{ij}\in\mathbf{Z}[1/p],n_2,n_3\in\mathbf{Z}.$$

This group is a variant of a construction by Abels of an interesting finitely presented solvable group, consisting of the $4\times 4$ southeast block of the above group (with $d=1$). Variations on Abels' group can also be found in \cite[Section 5.4]{CGP} and \cite{C1,Th}.

Let $M$, resp.~$M_\mathbf{Z}$, be the subgroup of $\Gamma$ consisting of matrices of the form

$$\begin{pmatrix}
  1 & 0 & 0 & 0 & m_1 \\
  0 & 1 & 0 & 0 & m_2 \\
  0 & 0 & 1 & 0 & 0 \\
  0 & 0 & 0 & 1 & 0 \\
  0 & 0 & 0 & 0 & 1 \\
\end{pmatrix},\quad m_1,m_2\in\mathbf{Z}[1/p],\text{ resp.}\in\mathbf{Z}.$$

The group $M\simeq\mathbf{Z}[1/p]^2$ is obviously normal in $\Gamma$, and the action of $\Gamma$ by conjugation on $M$ reduces to the action of $\textnormal{SL}_2(\mathbf{Z})$, so $M_\mathbf{Z}$ is normal as well in $\Gamma$. Theorem~\ref{main} follows from the following proposition, proved in Section \ref{pp}.

\begin{prop} The group $\Gamma/M_\mathbf{Z}$ satisfies the conditions of Theorem \ref{main}.\label{mainp}
\end{prop}

\medskip

\noindent {\bf Acknowledgment.} 
I thank Michel Coornaert who brought Gromov's question to my attention; I am grateful to Andreas Thom and Alain Valette who pointed out several mistakes in an early version of the paper.

%I thank Andreas Thom and Alain Valette for pointing out several mistakes in an early version of the paper.

\section{Elementary sofic groups}\label{es}

It is natural to introduce the class of elementary sofic groups as the smallest class of groups containing the trivial group $\{1\}$ and satisfying all four properties mentioned in the introduction. For instance, the group $\Gamma/M_\mathbf{Z}$ is elementary sofic (see Paragraph \ref{prosof}). We leave as open

\begin{que}
Is there any sofic group that is not elementary sofic?
\end{que}

Of course, it is very reasonable to conjecture a positive answer.

% For instance, the class of sofic groups is closed under free products \cite{ES} but I am unable to prove the same for the class of elementary sofic groups.

At least, we can give examples of non-(elementary sofic) groups. To this purpose, we have to give a constructive definition of this class. Let $\mathcal{C}_0$ be the class consisting of the trivial group only. By transfinite induction, define $\mathcal{C}_\alpha$ as the class of all groups in the following list
\begin{itemize}
\item groups all of whose finitely generated subgroups lie in $\bigcup_{\beta<\alpha}\mathcal{C}_\beta$;
\item finitely generated groups $G$ limits of a sequence $(G_n)$ with $G_n\in\bigcup_{\beta<\alpha}\mathcal{C}_\beta$ for all $n$;
\item groups $G$ in an exact sequence $1\to N\to G\to A\to 1$ with $A$ amenable and $N\in\bigcup_{\beta<\alpha}\mathcal{C}_\beta$.
\end{itemize}
Define $\mathcal{C}=\bigcup\mathcal{C}_\alpha$. Clearly, it consists of elementary sofic groups. To see the converse, we have to check

\begin{lem}
The class $\mathcal{C}$ is closed under taking subgroups and therefore exactly consists of elementary sofic groups.
\end{lem}
\begin{proof}
If $G$ lies in an extension of a group $N$ by an amenable group $A$, then any subgroup of $G$ lies in an extension of a subgroup of $N$ by a subgroup of $A$. The only other (straightforward) verification is that if $G$ is a limit of a sequence $(G_n)$ and $H$ is a finitely generated subgroup of $G$, then $H$ is a limit of subgroups of $G_n$. This shows at least that $\mathcal{C}$ is closed under taking finitely generated subgroups, but in turn this implies that $\mathcal{C}$ is closed under taking general subgroups. 
\end{proof}

\begin{lem}
Let $G$ be a non-amenable, finitely generated group. Assume that every finite index subgroup of $G$ is isolated, and that every amenable quotient group of $G$, is finite. Then $G$ is not elementary sofic.\label{coamn}
\end{lem}
\begin{proof}
Suppose that $G$ is elementary sofic. Let $N$ be a normal subgroup with $G/N$ amenable, with $N\in\mathcal{C}_\alpha$ with $\alpha$ minimal. As $N\in\mathcal{C}_\alpha$, and $\alpha\ge 1$ since $N$ is non-amenable, we have one of the three following possibilities
\begin{itemize}
\item We can write $N$ in an extension $1\to N'\to N\to A\to 1$ with $A$ amenable and $N'\in\bigcup_{\beta<\alpha}\mathcal{C}_\beta$. Since $N$ has finite index, the group $N'$ has finitely many conjugates only, so contains a subgroup $N''$ with $N''$ normal in $G$ and $G/N''$ is amenable as well. As $N''$ is a subgroup of $N'$, we have $N''\in\bigcup_{\beta<\alpha}\mathcal{C}_\beta$ as well, contradicting the minimality of $\alpha$. 
\item All finitely generated subgroups of $N$ belong to $\bigcup_{\beta<\alpha}\mathcal{C}_\beta$. Since $N$ is finitely generated, this would imply $N\in\bigcup_{\beta<\alpha}\mathcal{C}_\beta$ and cannot happen.
\item The group $N$ can be written as a non-trivial limit of groups. As $N$ is isolated by assumption, this cannot happen.\qedhere
\end{itemize}
\end{proof}

\begin{prop}
Let $G$ be a finitely presented group. Assume that $G$ is not residually finite, and that for some finite normal subgroup $Z$ of $G$, the group $G/Z$ is hereditary just infinite (every proper quotient of any finite index subgroup, is finite). Then $G$ is isolated, and if moreover $G$ is not amenable, then it is not elementary sofic.\label{nes}
\end{prop}
\begin{proof}
As $G$ is finitely presented, to prove that $G$ is isolated it is enough to prove that $G$ is finitely discriminable.
As $G$ is not residually 
finite, there exists an element $x\in G-\{1\}$ belonging to all finite index subgroups of $G$. We claim that $S=\{x\}\cup Z-\{1\}$ is a discriminating subset. Indeed, let $N$ be a normal subgroup of $G$ with $N\cap S=\emptyset$. If $N$ is contained in $Z$, then clearly $N=\{1\}$. Otherwise, the projection of $N$ into $G/Z$ is non-trivial, so has finite index, so $N$ itself has finite index. As $x\notin N$, this leads to a contradiction.

To prove that $G$ (now assumed non-amenable) is not elementary sofic, we note that the assumptions on $G$ are inherited by its finite index subgroups, and therefore the conditions of Lemma \ref{coamn} are satisfied. 
\end{proof}

Examples of groups satisfying the hypotheses of Proposition \ref{nes} are some non-residually finite lattices in (non-linear) simple Lie groups with finite center, see \cite[Paragraph 5.8]{CGP}, and non-amenable finitely presented simple groups, like Thompson's group $T$ (on the circle) or Burger-Mozes' groups \cite{BM}, which are amalgams of free groups. Accordingly, these groups are not elementary sofic; it is unknown if they are sofic.

The following proposition gives support that the class of elementary sofic groups is natural to introduce.

\begin{prop}                                                                    
The class of elementary sofic groups is closed under direct products and        
free products.                                                                  
\end{prop}                                                                      
\begin{proof}                                                                   
The case of direct products is straightforward and left to the reader.          
                                                                                
Let us prove that for every elementary sofic group $G$, the free product        
$G*F$ is elementary sofic for every free group $F$. It is enough to             
check that the class of                                                         
groups $G$ such that $G*F$ is elementary sofic contains $\{1\}$ (because        
free groups are residually finite, hence elementary sofic) and                  
satisfies the four stability properties mentioned in the introduction.          
The case of subgroups and direct limits is trivial. The case of limits          
is easy as well, since if $(G_n)$ tends to $G$, then $G_n*F$ tends to           
$G*F$ when $F$ is finitely generated; when $F$ is infinitely generated,         
this follows by a direct limit argument. It remains to prove that if $G$        
lies in an extension $$1\to N\to G\to M\to 1$$                                  
with $M$ amenable and $N*F$ is elementary sofic for every free group            
$F$, then                                                                       
$G*F$ is elementary sofic as well for every free group $F$. We                  
can write $G*F$ as $G\ltimes F^{*G}$ ($F^{*G}$ denoting the free product of copies of $F$ indexed by $G$). So the kernel of the natural map          
of $G*F$ onto $M$ is isomorphic to $N\ltimes F^{*G}$. If we fix a               
transversal of $N$ in $G/N$, we see that the latter group is isomorphic         
to $N\ltimes (F^{*G/N})^{*N}$, which in turn is isomorphic to                   
$N*(F^{*G/N})$. Since $F^{*G/N}$ is free, this group is elementary sofic by assumption, so $G*F$ is elementary          
sofic, by the fourth stability property.                                        
                                                                                
Now if $G$ and $H$ are elementary sofic groups, the free product $G*H$          
embeds into $(G\times H)*(G\times H)$, which in turn embeds into $(G\times H)*\mathbf{Z}$, since the latter is isomorphic to $\mathbf{Z}\ltimes         
(G\times H)^{*\mathbf{Z}}$. So $G*H$ is elementary sofic as well.                       
\end{proof} 

\section{Proof of Proposition \ref{mainp}}\label{pp}

Let $\Upsilon$ be the normal subgroup of $\Gamma$ consisting of elements for which $n_2=n_3=0$.
Define $\Lambda$ as the normal subgroup of $\Gamma$ consisting of elements for which 
$$\begin{pmatrix}
  a & b \\
  c & d \\
\end{pmatrix}=\begin{pmatrix}
  1 & 0 \\
  0 & 1 \\
\end{pmatrix}.$$

\subsection{The group $\Gamma/M_\mathbf{Z}$ is sofic}\label{prosof}

 As $\Gamma/\Upsilon\simeq\mathbf{Z}^2$, it is enough, to prove Proposition \ref{mainp}, to show that $\Upsilon/M_\mathbf{Z}$ is locally residually finite.

For $m\ge 0$, define $\Upsilon_m$ be the subset of $\Upsilon$ consisting of those matrices for which $$u_{02},u_{12},u_{23},u_{24}\in p^{-m}\mathbf{Z};$$ $$u_{03},u_{13},u_{24}\in p^{-2m}\mathbf{Z};$$ $$u_{04},u_{14}\in p^{-3m}\mathbf{Z}.$$
This is a subgroup, as can be check by a direct calculation, and $(\Upsilon_m)$ is clearly an increasing sequence of subgroups of union $\Upsilon$. Therefore it is enough for us to prove that $\Upsilon_m/M_\mathbf{Z}$ is residually finite.

Set $\Xi_m=\Upsilon_m\cap\Lambda$; this is a finitely generated nilpotent group, and $\Upsilon_m/M_\mathbf{Z}\simeq(\Xi_m/M_\mathbf{Z})\rtimes\textnormal{SL}_2(\mathbf{Z})$. Now it is known and easy to prove that any semidirect product of residually finite groups, where the normal factor is finitely generated, is also residually finite. So $\Upsilon_m/M_\mathbf{Z}$ is residually finite (the reader can check it is actually linear). 

\subsection{The group $\Gamma$ is isolated}

In \cite{CGP}, the following was proved.

\begin{prop}
A finitely generated group $G$ is isolated if and only if the two following conditions are fulfilled
\begin{itemize}
\item $G$ is finitely presented;
\item $G$ is finitely discriminable, i.e. there exists a finite ``discriminating" subset $X$ of $G-\{1\}$ such that every normal subgroup $N\neq\{1\}$ of $G$ satisfies $N\cap X\neq\emptyset$.
\end{itemize}
\end{prop}

Then $\Lambda/M_\mathbf{Z}$ is finitely discriminable. The proof is easy and strictly analogous to the case of the Abels group (case of the $4\times 4$ southeast block) \cite[Proposition~5.7]{CGP} and we skip it, just mentioning that we can pick $X$ as the set of elements of order $p$ in $M/M_\mathbf{Z}$.

Finite presentability of $\Lambda/M_\mathbf{Z}$ is a consequence of that of $\Lambda$, which is also analogous to the proof in the case of the Abels' group done in \cite{A1}, but requires slightly more work, so let us give the argument, not relying on the direct proof of \cite{A1} (that the reader can adapt), but on the general criterion of finite presentability from \cite{A2}.

Let $\Lambda(\mathbf{Q}_p)$ denote the same group as $\Lambda$, but with diagonal entries in $\mathbf{Q}_p^*$ and other entries in $\mathbf{Q}_p$. Let $\mathfrak{u}$ denote the Lie algebra of the unipotent part of $\Lambda(\mathbf{Q}_p)$. An element $(n_2,n_3)$ of $\mathbf{Z}^2$ is called a weight of $\mathfrak{u}$ if there exists $x\in\mathfrak{u}$ such that $\text{Ad}(d_i)(x)=\lambda_i x$ with $\log|\lambda|=n_i$ for $i=2,3$, where $d_2$, resp.~$d_3$, is the diagonal matrix $(1,1,p,1,1)$, resp. $(1,1,1,p,1)$. 

Abels proves that a necessary and sufficient criterion for finite presentability of $\Lambda$ is that the two following conditions are satisfied
\begin{itemize}
\item For any two weights $\alpha,\alpha'$ on $\mathfrak{u}/[\mathfrak{u},\mathfrak{u}]$, the segment in $\mathbf{R}^2$ joining $\alpha$ to $\alpha'$ does not contain 0;
\item Zero is not a weight on the second homology group $H_2(\mathfrak{u})$.
\end{itemize}

%These two conditions were checked in Lemmas 3.2 and 3.6 from \cite{C1}, where a strictly more general definition of $\mathfrak{u}$ was considered.

Let $(e_{ij})$ denotes the obvious basis of $\mathfrak{u}$ by elementary matrices.
A basis for the vector space $\mathfrak{u}/[\mathfrak{u},\mathfrak{u}]$ is given by the eigenvectors $e_{02}$, $e_{12}$, $e_{23}$, $e_{34}$ and we see immediately that the corresponding weights are $(1,0)$ (twice), $(-1,1)$, $(0,-1)$, so the first condition is satisfied.

The second homology group $H_2(\mathfrak{u})$ of $\mathfrak{u}$ is defined as
$\textnormal{Ker}(d_2)/\textnormal{Im}(d_3)$, where the maps
$$\mathfrak{u}\wedge\mathfrak{u}\wedge\mathfrak{u}\stackrel{d_3}
\to\mathfrak{u}\wedge\mathfrak{u}\stackrel{d_2}\to\mathfrak{u}$$ 
are defined by:

$$d_2(x_1\wedge x_2)=-[x_1,x_2]\;\;\text{and}$$ $$ d_3(x_1\wedge x_2\wedge
x_3)=x_3 \wedge [x_1,x_2]+x_2\wedge [x_3,x_1]+x_1\wedge
[x_2,x_3].$$

The weights of $\mathfrak{u}$ are the four ones described above, and the ones corresponding to the eigenvectors $e_{03}$, $e_{13}$, $e_{24}$, $e_{04}$, $e_{14}$ are $(0,1)$ (twice), $(-1,0)$, and $(0,0)$ (twice). If $e,e'$ are eigenvectors of weight $\alpha,\alpha'$, then $e\wedge e'$ has weight $\alpha+\alpha'$. Therefore the subspace of $\mathfrak{u}\wedge\mathfrak{u}$ corresponding to the weight 0 possesses as a basis the five elements
$$e_{i2}\wedge e_{24},\;e_{i3}\wedge e_{34}\;(i=0,1),\;e_{04}\wedge e_{14}.$$

It follows that the subspace of $\textnormal{Ker}(d_2)$ corresponding to the weight 0 possesses as a basis the three elements
$$e_{i2}\wedge e_{24}-e_{i3}\wedge e_{34}\;(i=0,1),\;e_{04}\wedge e_{14},$$
so to prove that $0$ is not a weight on $H_2(\mathfrak{u})$, we just have to check that these three elements belong to $\textnormal{Im}(d_3)$, for instance

$$e_{i2}\wedge e_{24}-e_{i3}\wedge e_{34}=d_3(e_{i2}\wedge e_{23}\wedge e_{34});$$
$$e_{04}\wedge e_{14}=d_3(e_{12}\wedge e_{24}\wedge e_{04}),$$
and the proof of Proposition \ref{mainp}, and therefore of Theorem \ref{main}, is complete.

% ----------------------------------------------------------------
\end{document}